\documentclass[11pt,reqno]{amsart}
\usepackage{amsmath,amsthm,amsfonts,amssymb,amscd,latexsym, calrsfs}
\usepackage[all]{xy}

\usepackage{amssymb}
\usepackage{mathrsfs}

\numberwithin{equation}{section}

\def\F{\mathbb{F}}

\def\su{{\subseteq}}
\def\la{{\langle}}
\def\ra{{\rangle}}

\def\w{{\omega}}

\def\dim{{\rm dim\,}}
\def\char{{\rm char}}

\def\Proof{\noindent{\sl Proof.}\ }
\def\qed{{\hfill $\Box$ \medbreak}}

\newtheorem{defi}{Definition}[section]
\newtheorem{thm}[defi]{Theorem}
\newtheorem{lem}[defi]{Lemma}
\newtheorem{cor}[defi]{Corollary}

\newtheorem{prop}[defi]{Proposition}
\newtheorem{rem}[defi]{Remark}

\usepackage[top=3.4cm,bottom=3.4cm,left=3.4cm,right=3.4cm, textwidth = 5in]{geometry}

\begin{document}

\title[Principal ideal enveloping algebras]{Enveloping algebras that are principal ideal rings}

\author{\textsc{Salvatore Siciliano}}
\address{Dipartimento di Matematica e Fisica ``Ennio De Giorgi", Universit\`{a} del Salento,
Via Provinciale Lecce--Arnesano, 73100--Lecce, Italy}
\email{salvatore.siciliano@unisalento.it}

\author{\textsc{Hamid Usefi}}
\address{Department of Mathematics and Statistics,
Memorial University of Newfoundland,
St. John's, NL,
Canada, 
A1C 5S7}
\email{usefi@mun.ca}

\thanks{The research of the second author was supported by NSERC of Canada under grant \# RGPIN 418201}

\begin{abstract}  Let $L$ be a restricted Lie algebra over a field of positive characteristic.
We prove that the restricted enveloping algebra of $L$ is a principal ideal ring if and only if
$L$ is an extension of a finite-dimensional torus by a cyclic restricted Lie algebra.
\end{abstract}

\date{\today}

  \maketitle

\section{Introduction}

Let $R$ be a ring with identity. 
Recall that   $R$  is called a \emph{principal right  ideal ring} (\emph{pri-ring} for short) if  every right ideal of $R$ is principal.  It is clear that such rings are right Noetherian and furthermore this property is inherited by homomorphic images. Similarly one defines  \emph{principal left ideal rings} (\emph{pli-rings}). Although there are examples of pri-rings that are not pli-rings and vice versa (see e.g. \cite{Lam}, \S 1, Example 1.25), the two properties turn out to be equivalent provided $R$ has an involution. For instance, this is the case when $R$ is a group algebra or an (ordinary or restricted) enveloping algebra. 
If a ring  is both a principal
right and left ideal ring then we simply call it a \emph{principal ideal ring}. 
It is known that a commutative principal   ideal ring is a finite direct sum of rings, which are either integral domains or are completely primary. 
Ore \cite{O} then was interested in certain differential polynomials which constitute a non-commutative principal ideal ring. This work was picked up and studied extensively by numerous algebraists including  Asano \cite{As}, 
Amitsur \cite{Am}, Jacobson \cite{J43}, and others. Most notably, Goldie in \cite{G} proved that 
a semiprime pri-ring is a finite direct sum of prime pri-rings and a prime pri-ring is a full matrix ring $M_n(K)$, where $K$ is a right noetherian integral domain. He further proved that
a left Noetherian pri-ring is a finite direct sum of pri-primary rings. Later Johnson in \cite{Jo} dropped the Noetherian assumption in Goldie's Theorem and proved that 
a ring is a pri-ring  if and only if it is a finite direct sum of primary pri-rings.
 
There has been also significant attention to group rings. First, Morita in \cite{Morita} characterized those finite groups whose group algebras over an algebraically closed field are principal ideal rings. Fisher and Sehgal in \cite{FS} extended Morita's result to 
nilpotent groups over any field, and later Passman \cite{P} dropped the nilpotence assumption from Fisher-Sehgal's result. We summarize the group ring result in the following. We denote the augmentation ideal of a group ring 
$K G$ by $\w(K G)$  and, for a prime $p$,  a finite group is called a $p^\prime$-group if its order is not divisible by $p$.

\begin{thm}[\cite{FS, P}]\label{FSP}
Let $K G$ be the group ring of $G$  over a field $K$. The following statements are equivalent.
\begin{enumerate}
\item $K G$ is a  principal right  ideal ring.
\item $K G$ is right Noetherian and   $\w(K G)$ is  principal as a right ideal.
\item   
\begin{itemize}
\item[] $\char\,  K = 0 $: $G$ is finite or finite-by-infinite cyclic.
\item[] $\char\, K = p> 0$:  $G$ is finite $p'$-by cyclic $p$ or 
finite $p'$-by infinite cyclic.
\end{itemize}
\end{enumerate}
\end{thm}

We just mention that among many other results,  Farkas and  Snider in \cite{FaSn} determined when the  augmentation ideal of the group ring $R G$ is principal as a right ideal, where $R$ is a commutative integral domain of characteristic 0, and a characterization of semigroup algebras that are pri-rings was obtained by Jespers and Okni\'nski in \cite{JO}.

In this paper we settle the same problem for another important  class of Hopf algebras.
Let $L$ be a restricted Lie algebra over a field $\F$ of characteristic $p>0$. The restricted universal enveloping algebra of $L$ is denoted by $u(L)$. We will characterize $L$ when $u(L)$ is a
principal ideal ring. Before stating our main result, we recall that an abelian restricted Lie algebra $T$ is called a torus if,  for every $x\in T$, the restricted subalgebra generated by $x^{[p]}$ contains $x$. Our main result is as follows.

\begin{thm}\label{us}
Let $L$ be a  restricted Lie algebra over a field of positive characteristic. Then  $u(L)$ is a principal ideal ring if and  only if  $L$ is an extension of a finite-dimensional  torus by a cyclic restricted Lie algebra.
\end{thm}

In other words, $u(L)$ is a principal ideal ring if and and only if  there exists a finite-dimensional torus $T$ such that $T$ is an ideal of $L$ and $L/T$ is cyclic as a restricted Lie algebra. 
In particular, we deduce from Theorem \ref{us} that if $u(L)$ is a principal ideal ring then $L$ is abelian. We note that, by Theorem \ref{FSP}, this is  not the case for group algebras, mainly due to the fact that semisimple group algebras are not necessarily commutative. 

We briefly explain the strategy of the proof. The sufficiency part is not difficult. In order to prove the necessity, we first deal with  the case that $L$ is finite-dimensional in Theorem \ref{fd} and consider the restricted ideal $T$ given by the intersection of $L$ with the last power of the augmentation ideal of $u(L)$. We show that $u(T)$ is  semisimple by proving that the counit of this Hopf algebra is not vanishing on the subspace of left integrals. By a well-known result of Hochschild [11], this allows us to conclude that $T$ is a torus. Then the proof of the  general case will use Theorem \ref{fd} in a crucial way.  

We  also show that the ordinary enveloping algebra $U(L)$ of an arbitrary Lie algebra $L$ is a principal ideal ring if and only if $L$ is either zero or 1-dimensional. In a more general setting, we think that finding the conditions under which an arbitrary Hopf algebra is a principal ideal ring would be an interesting future problem, and probably a difficult one.

\section{Preliminaries}
Let $L$ be a restricted Lie algebra over a field $\F$ of characteristic $p>0$. 
The terms of the lower  central series of $L$ are defined by  $\gamma_1(L)=L$ and 
$\gamma_{n+1}(L)=[\gamma_n(L), L]$, for every $n\geq 2$. We write $L^\prime$ for $\gamma_2(L)$ and $Z(L)$ for  the   center of $L$.   The centralizer of an element $x\in L$
is denoted by $C_L(x)$. For a subset $S$ of $L$, we denote by $\la S\ra_p$ the restricted subalgebra of $L$ generated by $S$ and by $\la S\ra_{\F}$ the subspace spanned by $S$.   Also, we denote by $S^{[p]^n}$ the restricted  subalgebra generated by all $x^{[p]^n}$, where $x\in S$. 
Recall that $S$ is called \emph{$p$-nilpotent}  if there exists an integer $n$ such that $S^{[p]^n}=0$.
An element $x\in L$ is called \textit{$p$-algebraic} if $\la x\ra_p$ is finite-dimensional and  
\textit{$p$-transcendental} otherwise. A restricted Lie algebra $L$ is said to be \emph{cyclic} if there exists an element $x$ of $L$ such that $L=\la x \ra_p$. If, furthermore, the element $x$ is $p$-nilpotent then $L$ is  called \emph{nilcyclic}. Following \cite{Hoch}, we say that a restricted Lie
algebra is \emph{strongly abelian} if it is abelian and its power mapping is zero.

We  denote by $u(L)$  the restricted (universal) enveloping algebra of $L$.
A restricted Lie algebra $T$ is said to be a {\it torus} if  $T$ is abelian and every element of $T$ is semisimple. It is  well known, by a  theorem of 
Hochschild \cite{H},  that $u(L)$  is a semisimple algebra   if and only if $L$ is a finite-dimensional  torus. 
The Poincar\'e-Birkhoff-Witt  (PBW) Theorem  for restricted Lie algebras 
(see e.g. \cite[Chapter 2, \S 5, Theorem 5.1]{SF}) states that if $X$ is an ordered basis of $L$ then $u(L)$ has a basis consisting of 
PBW monomials of the form
$$
x_1^{a_1}\cdots x_n^{a_n},
$$
where the $x_1\leq \cdots \leq x_n$ in $X$ and $0\leq a_i\leq p-1$, for all $i$. Let $H$ be a restricted subalgebra of $L$. 
Recall  that $u(L)$ is a free left $u(H)$-module. Indeed, if  $Y$ is  a basis  of a complementary vector subspace of $H$ in $L$, then 
\begin{align*}
u(L)= \oplus u(H) w,
\end{align*} 
where the sum is over all the PBW monomials $w$ in $Y$.

The following result from \cite{P} will be used in the sequel  and we quote it here 
for  convenience of the reader.

\begin{lem}\label{passman}
Let $R$ be a prime ring which is both right and left
Noetherian. Let $A\neq  R$ be a 2-sided ideal with $A = \alpha R$ for some $\alpha \in R$
and let $B$ be another 2-sided ideal of $R$. Then
\begin{enumerate}
\item[i.] $AB = B$ implies that $B = 0$.
\item[ii.] If $B$ is prime and $B \subsetneq  A$, then $B = 0$.
\end{enumerate}
\end{lem}

For a subset $S$ of a ring $R$,  we  denote by $ \mathrm{ann}_R^r(S)$  the right  annihilator of $S$ in $R$.
Finally, we recall that an element $t$  of a Hopf algebra $H$ is said to be a \emph{left} (respectively, \emph{right}) \emph{integral} if $h t=\epsilon (h)t $ (respectively, $t h=\epsilon (h)t$), for every $h \in H$, where $\epsilon$ denotes the counit of $H$. We will denote by $\int_{H}^l$ and $\int_{H}^r$, respectively,  the subspaces of left and right integrals of $H$. Clearly, $\int_{H}^l$  is a left ideal and $\int_{H}^r$ is a right ideal of $H$ and, when $H$ is finite-dimensional, they are both 1-dimensional. We shall also denote by $\w(H)$ the augmentation ideal of $H$, that is, the kernel of $\epsilon$.  In case $H=u(L)$ we
write  $\w(L)$ for $\w(u(L))$. 
Note  that $\w(L)=Lu(L)=u(L)L$. For every positive integer $n$, we will consider the restricted ideal $D_{n}(L)$ of $L$ defined in  
\cite{RS2} as follows:
$$
D_{n}(L)=L\cap \w(L)^{n}=\sum_{ip^{j}\geq n}{{\gamma}_{i}(L)}^{[p]^{j}},
$$
where $ \w(L)^{n}$ is the $n$-th power  of the ideal $\w(L)$.

 For more notation and some well-known results about restricted Lie algebras and Hopf algebras we refer the reader to \cite{SF} and \cite{DNR, M}, respectively.

\section{Proof of the main result and concluding remarks}

Throughout, all restricted Lie algebras are defined over a field $\F$ of positive  characteristic $p$.

\begin{lem}\label{ann-hopf}
		Let $\mathcal{H}$ be a Hopf algebra over a field  and $\mathcal{B}$ a finite-dimensional 
		Hopf subalgebra of $\mathcal{H}$. Suppose further that  $\mathcal{H}$ is a free left  
		$\mathcal{B}$-module.  Then   
		$\mathrm{ann}_{\mathcal{H}}^r(\w(\mathcal{B}))=\int_{\mathcal{B}}^l  \mathcal{H}$.
\end{lem}
\Proof We clearly have  
	$\int_{\mathcal{B}}^l  \mathcal{H}\su \mathrm{ann}_{\mathcal{H}}^r(\w(\mathcal{B}))$. Since   $\mathcal{H}$ is a free left $\mathcal{B}$-module, we have
\begin{align}\label{free}
\mathcal{H}= \oplus \mathcal{B} w_i,
\end{align} 
where the  $w_i$'s form a $\mathcal{B}$-basis of $\mathcal{H}$. Let $v\in\mathrm{ann}_{\mathcal{H}}^r(\w(\mathcal{B}))$.
Then, by Equation \eqref{free}, we have $v=\sum u_i w_i$, where each $u_i$ is in $\mathcal{B}$.
Now, we have 
$$
xv=\sum xu_i w_i=0,
$$
for every $x\in \mathcal{B}$. We deduce from Equation \eqref{free} that $ xu_i=0$, for every $x\in \mathcal{B}$ and all the $u_i$'s.
Consequently, each $u_i$ annihilates $\w(\mathcal{B})$. Since $\dim \int_{\mathcal{B}}^l=1$, it follows that   the right annihilator of 
$\w(\mathcal{B})$ in $\mathcal{B}$ is just  $\int_{\mathcal{B}}^l$. Thus, we have  $v=\sum u_i w_i\in \int_{\mathcal{B}}^l  \mathcal{H}$, as required.
\qed

Let $\mathcal{H}$ be a Hopf algebra over a field  and $\mathcal{B}$ a finite-dimensional 
Hopf subalgebra of $\mathcal{H}$. Note that Lemma \ref{ann-hopf} applies in particular when    $\mathcal{H}$ is 
finite-dimensional (see \cite{NZ})  or pointed (see \cite{Rad}), or $\mathcal{B}$ is 
semisimple (see \cite{NZ92}). Since restricted enveloping algebras are pointed, we deduce the following:

\begin{cor}\label{ann}
Let $L$ be a  restricted Lie algebra and $H$ a finite-dimensional restricted 
subalgebra of $L$. Then 
$\mathrm{ann}_{u(L)}^r(\w(H))=\int_{u(H)}^l  u(L)$.
\end{cor}

\begin{lem}\label{torus}
Let $L$ be a finitely generated abelian restricted Lie algebra. If $L= L^{[p]}$ then $L$ is a finite-dimensional torus.
\end{lem}
\Proof Let $\mathcal{L}=L\otimes_{\F}\bar{\F}$, where $\bar{\F}$ denotes the algebraic closure of  $\F$. Note that  $\mathcal{L}$ is $\bar \F$-spanned by elements of the form $x^ {[p]}\otimes 1$ with $x\in L$, hence we still have  $\mathcal{L}= \mathcal{L}^{[p]}$.
Also,  $\mathcal{L}$ is   finitely generated.   We claim that $\mathcal{L}$ is finite-dimensional. Indeed, by well known results about the structure of finitely generated abelian restricted Lie algebras over perfect fields (see e.g. Section 4.3 in \cite{BMPZ}), we deduce that there exist  $x_1, \ldots, x_r\in \mathcal{L}$ such that  
$$
\mathcal{L}= \la  x_1 \ra_p \oplus \cdots \oplus \la  x_r \ra_p.
$$ 
Now, each $x_i$ is in $\mathcal{L}^{[p]}$ and by the decomposition above we have that 
$x_i\in  \la  x_i^{[p]} \ra_p$. This means that all the $x_i$'s are $p$-algebraic and hence 
$\mathcal{L}$ is finite-dimensional. Thus, $L$ is finite-dimensional and it follows from  
 \cite[\S 2.3, Proposition 3.3]{SF} or  \cite[Proposition 4.5.4]{W}
that $L$ is a torus. 
\qed

We are now ready to settle the finite-dimensional case.

\begin{thm}\label{fd}
Let $L$ be a finite-dimensional restricted Lie algebra. If $u(L)$ is a pri-ring then $L$ is an extension of a torus by a nilcyclic restricted Lie algebra. In particular, $L$ is abelian.
\end{thm}
\Proof
As $u(L)$ is a pri-ring,  there exists $w\in \w(L)$ such that $\w(L)=w u(L)$. Consider the onto map $\rho: u(L)\to u(L) w$ given by   right multiplication by $w$.  We claim that 
$\ker \rho$ has codimension 1 in $u(L)$. To prove the claim, we recall  that $\dim \int_{u(L)}^r=1$
and  $\int_{u(L)}^r\su \ker\rho$. On the other hand, let $z\in \ker\rho$. Then we have 
$z \w(L)=0$. Hence, for every $u\in u(L)$, we get $zu=\epsilon (u) z$. We deduce that  $z \in  \int_{u(L)}^r$. Hence, $\ker \rho=\int_{u(L)}^r$ and this yields the claim. Now it follows that 
$$
\dim \w(L)=\dim u(L) -1=\dim u(L) w.
$$
Since  $u(L) w\su \w(L)$, we deduce from the above that $w u(L)=\w(L)= u(L) w$.
Hence, for every positive integer $n$, we have   $\w^n(L)= u(L) w^n=  w^n u(L) $. 
Since $L$ is finite-dimensional, there exists $n$ such that $\w^n(L)= \w^{n+1}(L)$ and, in turn, 
there exists an element $u\in u(L)$ such that $w^n= w^{n+1} u $. Hence, 
 $w^n(1-wu)=0$ and subsequently $\w^n(L)(1-wu)=0$. Let $T=D_n(L)$. Then we have 
 $\w(T)(1-wu)=0$. Let $t$ be a non-zero left integral of the Hopf algebra $u(T)$.
We can see, by Lemma \ref{ann},  that $1-wu=t v$, for some $v\in u(L)$. Since $w\in \w(L)$ we deduce that
$\epsilon(t)\neq 0$. It now follows from \cite{LS} that $u(T)$ is semisimple.
So,  by Hochschild's Theorem \cite{H}, the restricted ideal $T$ is a torus.
Note that $L^{[p]^{j}}\su T$ whenever $p^j\geq n$. Hence,  $L/T$ is $p$-nilpotent. Now, we claim  that $H=L/T$ is  cyclic.
We observe that the Frattini subalgebra $\Phi_p(H)$ of $H$ is  equal to $D_2(H)=H'+H^{[p]}$ (see \cite[Corollary 5.2(ii)]{LT}). Let $\bar H=H/ D_2(H)$. Note that $u(\bar H)$ is  a pri-ring. Moreover, $\bar{H}$ is strongly abelian and  
 $u(\bar H)$ is isomorphic to the truncated polynomial ring  
$\F[X_1, \ldots, X_r]/\la X_1^p, \ldots, X_r^p\ra$, where $r=\dim \bar H$. But it is easy to see that the truncated polynomial ring  is  a principal ideal ring if and only if $r=1$. Therefore, by \cite{LT}, we deduce that $H$ is cyclic which proves the claim. Since $T$ is a torus, we have that $z\in \la z^{[p]}\ra_p$, for every $z\in T$. As $T$ is also an abelian  restricted ideal,   it follows that $T$ is central in $L$.  Finally, since $L/T$ is cyclic, we conclude  that $L$ is abelian.\qed

\begin{lem}\label{fg}
Let $L$ be a  restricted Lie algebra. If $u(L)$ is (right) Noetherian then    every restricted subalgebra of $L$ is finitely generated.
\end{lem}
\Proof  Let $H$ be a restricted subalgebra of $L$.  Suppose to the contrary that  $H$ is not finitely generated. Then we can find  $x_1, x_2,\ldots \in H$ so that each $x_{k+1}$ is not contained in $\la  x_1, \ldots,  x_k \ra_p$. For every $k\geq 1$, let $H_k=\la x_1, x_2 \cdots x_k \ra_p$. 
Now consider the following  ascending   chain of right ideals of $u(L)$:
\begin{align*}
H_1 u(L)\subseteq H_2 u(L)\subseteq \cdots \subseteq 
 H_k u(L)\subseteq\cdots .
\end{align*}
Then there exists an integer
$n$ such that $H_n u(L)=H_{n+1} u(L)$. It follows that $x_{n+1}\in L\cap H_n u(L)$. Since
$u(L)$ is a free left $u(H_n)$-module, we have that 
$L\cap H_n u(L)=H_n$ and so $x_{n+1}\in H_n$,  a contradiction.\qed

It is still unkwown when $u(L)$ is a domain \cite[Problem 3.59]{Dniester}. More precisely, if $u(L)$ is a domain then it is clear that $L$ has no nonzero $p$-algebraic elements, however it is an open problem if the converse is also true. In the following result, we establish when $u(L)$ is a principal ideal commutative domain.

\begin{prop}\label{ab}
Let $L\neq 0$ be an abelian restricted Lie algebra. Then $u(L)$ is a principal ideal domain if and only if $L$ is infinite dimensional cyclic.
\end{prop}
\Proof If $L$ is infinite dimensional cyclic, then $u(L)$ is isomorphic to a polynomial algebra in one indeterminate over $\F$ and so it is a principal ideal domain. Conversely, suppose that $u(L)$ is a principal ideal domain. It is clear that every nonzero element of $L$ is $p$-transcendental. Let $T=\cap_{i=1}^{\infty} L^{[p]^i}$. Since every pri-ring  is right Noetherian,  we deduce from Lemma \ref{fg} that $T$ is finitely generated. Consequently, Lemma  \ref{torus} implies that $T$ is a finite-dimensional torus and so $T=0$. Denote by $\bar{\F}$ the algebraic closure of $\F$ and let $\bar{L}=L\otimes_{\F}\bar{\F}$. Since, by Lemma \ref{fg},  $L$ is finitely generated, so is $\bar{L}$.  Therefore, by Theorem 3.1 of \cite[Chapter 4, \S 3]{BMPZ}, we see that 
$$
\bar{L}=\langle \bar x_1\rangle_p \oplus \cdots \oplus \langle \bar x_r \rangle_p,
$$
 for suitable $\bar x_1,\ldots,\bar x_r\in \bar L$. As $L$ is infinite dimensional, at least one of such elements $\bar x_i$ is transcendental, say $\bar x_1$. Now, suppose that $\bar x_j$ is $p$-algebraic. Then,  by Theorem 4.5.8 of \cite{W} we have that $\langle \bar x_j \rangle_p=\mathfrak{T}\oplus \mathfrak{N}$, where  $\mathfrak{T}$ is a torus and 
$\mathfrak{N}$ a $p$-nilpotent restricted Lie algebra. Since
 $$
 \cap_{i=1}^{\infty} \bar{L}^{[p]^i}=\cap_{i=1}^{\infty}( L^{[p]^i}\otimes_{\F}\bar{\F})=\left( \cap_{i=1}^{\infty} L^{[p]^i}\right) \otimes_{\F}\bar{\F}=0
 $$
  and  $\mathfrak{T}^{[p]}= \mathfrak{T}$, we must have  $\mathfrak{T}=0$, in particular $\bar x_j$ is $p$-nilpotent. As a consequence, we have that $r=\dim_{\bar{\F}}\bar{L}/\bar{L}^{[p]}\leq \dim_{\F}L/L^{[p]}$.  On the other hand, as $u(L/L^{[p]})$ is   a principal ideal ring, it follows from Theorem \ref{fd}  that $L/L^{[p]}$ is nilcyclic. This forces $r=1$, so that $\bar{L}=\langle \bar x_1 \rangle_p$ is infinite dimensional cyclic.   Now, as $u(L)$ is a principal ideal domain, the ideal $\omega(L)$ is generated by an element $v\in u(L)$. It follows that $v\otimes 1$ generates $\omega(\bar{L})$ as an ideal of $u(\bar{L})$ and so $v\otimes 1$  is proportional to $\bar x_1$.  In particular, $\bar{L}$ is spanned by all the $p$-powers of $v\otimes 1$. This allows us to conclude that $v\in L$ and  $L=\langle v\rangle_p$, completing the proof.
 \qed

In our next result, we  deal with the case that  $L$ is abelian of infinite dimension.
 
\begin{prop}\label{ab-inf}
Let $L$ be  an infinite dimensional abelian restricted Lie algebra. If $u(L)$ is a principal ideal ring then  $L$ is an extension of a torus by an infinite dimensional cyclic restricted Lie algebra.
\end{prop}
\Proof Let $T=\cap_{i=1}^{\infty} L^{[p]^i}$ and set $\mathcal{L}=L/T$. We know by Lemma \ref{fg} that both  $L$ and $T$ are  finitely generated.
Now  Lemma \ref{torus} implies  that $T$ is a finite-dimensional torus. Hence the restricted Lie algebra $\mathcal{L}$ is infinite dimensional. Now, if $\mathcal{L}=\mathcal{L}^{[p]}$ then we infer from Lemma \ref{torus} that $\mathcal{L}$ is  finite-dimensional, which is a contradiction. Therefore we have  $\mathcal{L}^{[p]}\subsetneq \mathcal{L}$. 
Moreover, by Theorem \ref{fd}, each quotient $\mathcal{L}/ \mathcal{L}^{[p]^i}$ is nilcyclic and so 
$u(\mathcal{L}/ \mathcal{L}^{[p]^i})\cong \F[X]/(X^{p^i})$. For every  $i\leq j$, we have the natural
projection $\mathcal{L}/ \mathcal{L}^{[p]^{j}}\to \mathcal{L}/ \mathcal{L}^{[p]^i}$ which induces 
the projections 
\begin{align*}
u(\mathcal{L}/ \mathcal{L}^{[p]^{j}})\to u(\mathcal{L}/ \mathcal{L}^{[p]^i}).
\end{align*}
We observe that $\varprojlim u(\mathcal{L}/ \mathcal{L}^{[p]^i})\cong \F[[ X]]$, where  $\F[[ X]]$ is the  formal power series over $\F$ (see e.g. \cite[Chapter III, \S 10]{Lang}).   Furthermore, as $\cap_{i=1}^{\infty} \mathcal{L}^{[p]^i}=0$ we have 
the natural algebra embeddings  
\begin{align*}
u(\mathcal{L})\hookrightarrow u(\varprojlim \mathcal{L}/ \mathcal{L}^{[p]^i})\hookrightarrow \varprojlim u(\mathcal{L}/ \mathcal{L}^{[p]^i})\cong \F[[ X]].
\end{align*}
Since $\F[[ X]]$ is a domain, we deduce that $u(\mathcal{L})$ is a principal ideal domain. Thus, Proposition \ref{ab} implies that  $\mathcal{L}=L/T$ is infinite dimensional cyclic and the proof is complete. \qed

 Following \cite{BP}, for a restricted Lie algebra $L$ we define
$$
\Delta(L)=\{x\in L \vert\, \dim [L,x]< \infty\}.
$$
Then $\Delta(L)$ is clearly a restricted ideal of $L$ which is the Lie algebra analogue of the FC-center of a group.

\begin{prop}\label{delta=0}
Let $L$ be a  restricted Lie algebra such that $\Delta(L)=0$. If $u(L)$ is a pri-ring then $L=0$.
\end{prop}
\Proof Suppose to the contrary that $L\neq 0$. It follows from Corollary 6.4 of \cite {BP} that $u(L)$ is a prime ring. 
By Lemma \ref{passman} we have that $ \w^{n+1}(L) \subsetneq \w^n(L)$, for every positive 
integer $n$. Let $N=\cap_{i=1}^{\infty} D_i(L)$. Since $\w(N)\su \w^n(L)$ for every $n$, 
it follows that $u(L/N)\cong u(L)/Nu(L)$ is infinite-dimensional.
Hence, $L/N$ is infinite-dimensional.  Note that, by Lemma \ref{fg},  $L$ is finitely generated and, consequently, $L/D_n(L)$ is finite-dimensional and $p$-nilpotent, for every $n$.
Therefore, by Theorem \ref{fd}, $L/D_n(L)$ is cyclic. In particular, we have  $L'\su D_n(L)$, for 
every $n$. Hence,  $L'\su N$ and $\bar L=L/N$ is abelian. Now, by Lemma \ref{ab-inf}, there exists a restricted ideal $T$ of $L$ containing $N$ such that $\bar T=T/N$ is a torus and $\bar L/\bar T\cong  L/ T $ is infinite dimensional cyclic. Let $x\in T$. Then, for every positive integer $n$, we have 
 
\begin{align*}
x=\sum \alpha_i y_i^{[p]^n}  \text{ modulo } N,
\end{align*}
 where each $y_i\in T$ and $\alpha_i\in \F$. Hence, 
$x\in D_m(L)$, for every $m$. We conclude that $T=N$ and so $ L/ N $ is infinite dimensional cyclic. As a consequence, we have that $\F[X] \cong u(L/N)\cong u(L)/Nu(L)$ is a domain. Note that $Nu(L)\subsetneq \w(L)$ and it follows from  Lemma  \ref{passman} that $N=0$. We conclude that $L$ is abelian and so  
$L=\Delta(L)$, a contradiction. \qed

The next result shows that our problem boils down to the commutative case.

\begin{prop}\label{infd}
Let $L$ be a  restricted Lie algebra. If $u(L)$ is a pri-ring then $L$ is  abelian.
\end{prop}
\Proof Let $\Delta_0=0$ and inductively define $\Delta_{n+1}$ by
$$
\Delta_{n+1}/\Delta_{n}=\Delta (L/\Delta_{n}), 
$$
for every $n\geq 0$. 
Then we obtain a chain $0\subseteq \Delta_1\subseteq \cdots\subseteq \Delta_n\subseteq \cdots$ of restricted ideals of $L$ which must stabilize, that is, there exists a smallest  positive integer $n$ such that  $\Delta_n=\Delta_{n+1}$. Since $u(L/\Delta_n)$ is  a pri-ring and 
$\Delta(L/\Delta_n)=0$,   Proposition \ref{delta=0} implies that $L=\Delta_n$.
Suppose first that $n=1$, that is, $L=\Delta(L)$. By Lemma \ref{fg},  there exist $x_1, \ldots, x_r\in L$ such that $L=\la x_1, \ldots, x_r\ra_p$. We have 
$$
Z(L)=\cap_{i=1}^r C_L(x_i).
$$
 Note that 
$\dim L/C_L(x_i) =\dim [L, x_i]$ is finite-dimensional, for every $1\leq i\leq r$. Therefore  
$L/Z(L)$ is also finite-dimensional. We deduce from Theorem \ref{fd} that $H=L/Z(L)$ is abelian and so $L$ is nilpotent of class at most 2. It follows that $L^{[p]}\su Z(L)$. Thus, 
$u(H)$ is a truncated polynomial ring and since $u(H)$ is also a pri-ring, we deduce that $\dim H=1$.
It now follows that $L$ is abelian. 

Now suppose, if possible, that $n\geq 2$ and let 
$\mathcal{L}= L/\Delta_{n-2}$. Note that $\Delta_2(\mathcal{L})=\mathcal{L}$. It follows from what we proved above that $\mathcal{L}/\Delta(\mathcal{L})$ is abelian. By Lemma \ref{fg}, there 
exist $y_1, \ldots, y_r\in \mathcal{L}$ such that 
$$
\Delta(\mathcal{L})=\la y_1, \ldots, y_r\ra_p.
$$
Notice that 
$$
[\mathcal{L}, \Delta(\mathcal{L})]=\sum_{i=1}^r [\mathcal{L}, y_i].
$$
Since each subspace $[\mathcal{L}, y_i]$ is finite-dimensional, we deduce that
$[\mathcal{L}, \Delta(\mathcal{L})]$ is also finite-dimensional. 
Let $x\in \mathcal{L}$. Note that 
$$
[\mathcal{L}, x^{[p]}]\su [\cdots [\mathcal{L}, \underbrace{x], \ldots, x}_p]\su [\Delta(\mathcal{L}), x]\su [ \Delta(\mathcal{L}), \mathcal{L}]
$$ is finite-dimensional. Hence, $x^{[p]}\in \Delta(\mathcal{L})$.
It follows that  $\mathcal{L}^{[p]}\su \Delta(\mathcal{L})$ and $\mathcal{L}/\Delta(\mathcal{L})$ is strongly abelian.
Since 
$$
u(\mathcal{L}/\Delta(\mathcal{L}))\cong u(\mathcal{L})/\Delta(\mathcal{L})u(L)
$$ is  a principal ideal ring, we deduce that 
$\mathcal{L}/\Delta(\mathcal{L})$ is at most 1-dimensional.
But now it easily follows that $[\mathcal{L}, \mathcal{L}]$ is finite-dimensional and so $\mathcal{L}=\Delta(\mathcal{L})$, a contradiction.
\qed

\smallskip
It is now a simple matter to prove the main result of the paper:

\smallskip
{\it Proof of Theorem \ref{us}:}
Suppose first that $u(L)$ is a principal ideal ring. Then, by Proposition \ref{infd}, $L$ is abelian.
The necessity part is then a combination of Theorem \ref{fd} and  Proposition \ref{ab-inf}. 
Conversely, suppose that $L$ contains a finite-dimensional toral ideal $T$ such that
$L/T$ is a cyclic restricted Lie algebra.
Then there exists $x\in L$ such that $L/T$ is generated as a restricted Lie algebra by
  $\bar  x$, where $\bar  x$ is the image of $x$ in $L/T$. Thus we have $L=T+ \la x\ra_p$.
 Since  $T$ is an abelian restricted ideal consisting of semisimple elements,  it follows that  $L$ is abelian. Indeed, let   $y\in T$. Then $y\in \la y^{[p]}\ra_p$. Hence,
 $[x,y]=\sum_{i\geq 1} \alpha_i [x,y^{[p]^i}]=0$.
  Now, suppose first that $\bar x$ is $p$-transcendental. Then $L=T\oplus \la x\ra_p$. It follows that
	\begin{equation}\label{tensor}
	u(L)\cong u(T)\otimes_{\F} u(\la x\ra_p)\cong u(T)\otimes_{\F} \F[X].
	\end{equation}
	Since, by Hochschild's Theorem \cite{H}, $u(T)$ is a commutative semisimple algebra, we
	have that  
	$$
	u(T)\cong \F_1 \oplus \cdots \oplus \F_r,
	$$
	where $\F_1, \ldots,\F_r$ are field extensions of the ground field $\F$. 
	Hence, from (\ref{tensor}), it follows that
	\begin{equation}\label{extensions}
	u(L)\cong \F_1[X] \oplus \cdots \oplus \F_r[X].
\end{equation}
	Since a finite direct sum of polynomial algebras is a principal ideal ring, we then conclude from (\ref{extensions}) that $u(L)$ is also a principal ideal ring.  Finally, suppose that $\bar x$ is $p$-algebraic. Let $\mathcal{L}$ be the direct sum of $T$ and   an infinite dimensional cyclic restricted Lie algebra $H$. As $H$ is a free restricted Lie algebra, it is clear that there exists a restricted homomorphism $\pi$ from $\mathcal{L}$ onto $L$. Thus 
$$
u(L)\cong u(\mathcal{L}/\ker \pi)\cong u(\mathcal{L})/(\ker \pi)u(\mathcal{L}).
$$
But we already proved that  $u(\mathcal{L})$ is  a principal ideal ring, so $u(L)$ is also a principal ideal ring, completing the proof.
\qed

\begin{rem}\emph{
Let $L$ be a restricted Lie algebra over a field of characteristic $p>0$ and suppose that $u(L)$ is a pri-ring. Then, by our Theorem \ref{us}, $L$ is an extension of a finite-dimensional  torus by a cyclic restricted Lie algebra. If $L$ is infinite dimensional, we showed in the proof  Theorem \ref{us} that this extension is indeed split. Furthermore, the same conclusion follows from \cite[Theorem 4.5.8]{W}  if   the ground field  is perfect and $L$ is finite-dimensional.   However, this might  not be the case in general.  For an explicit example, let $\F$ be a field of positive characteristic $p$ containing an element $\alpha$
with no $p$-th root in $\F$ and consider the abelian restricted Lie algebra $L=\F x+ \F y$ with $x^{[p]}=x$ and $y^{[p]}=\alpha x$. 
}
\end{rem}

Finally, let $L$ be a Lie algebra over an arbitrary field and denote by $U(L)$  its ordinary enveloping algebra. Suppose that $U(L)$ is a principal ideal ring and let $\Omega(L)$ be the augmentation ideal of $U(L)$. As $U(L)$ is a prime ring, if $\Omega^2(L)=\Omega(L)$ then Lemma \ref{passman} forces $L=0$. On the other hand, if $\Omega^2(L)\subsetneq \Omega(L)$ then  $L^\prime \subsetneq L$ and $U(L/L^\prime)\cong U(L)/L^\prime U(L)$ is a principal ideal ring which is  isomorphic to the polynomial algebra in $\dim L/L^{\prime}$ indeterminates. Hence $\dim L/L^{\prime}=1$ and, moreover, as $I=L^\prime U(L)=U(L)L^\prime$ is a 2-sided prime ideal of $U(L)$ with $I\subsetneq \Omega(L)$,  by Lemma \ref{passman} we conclude that $I=0$. Thus $L^\prime =0$ and $L$ is 1-dimensional. We thereby proved the following

\begin{thm} Let $L$ be a Lie algebra over any field. Then $U(L)$ is a principal ideal ring if and only if $L$ is either zero or 1-dimensional.
\end{thm}

\section*{Acknowledgements}
We are grateful to the referee for careful reading of the
manuscript and helpful comments. The second author would like to thank the
Dipartimento di Matematica e Fisica ``Ennio De Giorgi" of the Universit\`{a} del Salento  for its hospitality during
his visit while this work was completed.

\end{document}